\documentclass[letterpaper, 10 pt,twocolumn,twoside,journal]{IEEEtran}  

\IEEEoverridecommandlockouts                              

\usepackage{graphicx} 
\usepackage{epsfig} 
\usepackage{cite} 
\usepackage{amsmath} 
\usepackage{amssymb}  
\usepackage{mathrsfs}
\usepackage{mathtools}
\usepackage{color}
\usepackage{nicefrac}
\usepackage[linesnumbered, ruled, vlined]{algorithm2e}
\mathtoolsset{showonlyrefs}
\usepackage{amsthm}
\graphicspath{{.}}

\theoremstyle{plain}

\usepackage{etoolbox}
\usepackage{float}
\usepackage{tabularx}
\usepackage{enumitem}

\usepackage[font=small,skip=0pt]{caption}

\usepackage{color}

\theoremstyle{definition}
\newtheorem{dfn}{Definition}

\newtheorem{lemma}{Lemma}

\newtheorem{problem}{Problem}
\newtheorem{theorem}{Theorem}

\newtheorem{assumption}{Assumption}

\newcommand{\R}{\mathbb{R}}
\newcommand{\N}{\mathbb{N}}

\newcommand{\barsig}{\overline{\Sigma}}

\newcommand{\diag}{\textnormal{diag}}

\newcommand{\tr}{\textnormal{tr}}

\DeclareMathOperator*{\argmax}{arg\,max}
\DeclareMathOperator*{\argmin}{arg\,min}

\makeatletter
\patchcmd{\@makecaption}
{\scshape}
{}
{}
{}
\makeatother

\title{\LARGE \bf
Differentially Private LQ Control
}

\author{Kasra Yazdani$^{\ast}$, Austin Jones$^{\dagger}$, Kevin Leahy$^{\dagger}$, Matthew Hale$^{\ast}$
\thanks{$^{\ast}$Kasra Yazdani and Matthew Hale are with the Department of  Mechanical and Aerospace Engineering at the University of Florida, Gainesville, FL USA. Emails: \texttt{\{kasra.yazdani,matthewhale\}@ufl.edu.} 
	
$^{\dagger}$Austin Jones and Kevin Leahy are with the Massachusetts Institute of Technology Lincoln Laboratory, Lexington, MA USA. Emails: \texttt{\{austin.jones,kevin.leahy\}@ll.mit.edu.}} 

\thanks{KY and MH were supported in part by NSF CAREER Grant No. 1943275 and
by AFOSR under Grant FA9550-19-1-0169.}
\thanks{DISTRIBUTION STATEMENT A:Approved for public release. Distribution is unlimited.}%
\thanks{This material is based upon work supported by the United States Air Force under Air Force Contract No. FA8702-15-D-0001. Any opinions, findings, conclusions or recommendations expressed in this material are those of the author(s) and do not necessarily reflect the views of the United States Air Force.}%
}%

\begin{document}

\maketitle

\begin{abstract}
As multi-agent systems proliferate and share more user data, new approaches are needed to protect sensitive data while still enabling system operation. 
To address this need, this paper presents a private multi-agent LQ control framework. 
Agents' state trajectories can be sensitive and we therefore protect them using differential privacy. 
We quantify the impact of privacy along three dimensions: the amount of information shared under privacy, the control-theoretic cost of privacy, and the tradeoffs between privacy and performance. These analyses are done in conventional control-theoretic terms, which we use to develop guidelines for calibrating privacy as a function of system parameters. Numerical results indicate that system performance remains within desirable ranges, even under strict
privacy requirements. 
\end{abstract}
\section{Introduction}

\IEEEPARstart{M}{ulti-agent} systems, such as smart power
grids and robotic swarms, require agents
to exchange information to work together. In some cases,
the information shared may be sensitive. For example, consumption data in a power grid can expose habits and activities of individuals~\cite{McDaniel2000,doe2010}. 
Sensitive user data must be protected when it is shared, though of course it must remain useful in multi-agent coordination. Hence, privacy in  multi-agent control should protect sensitive data from 
the agent receiving it
while still ensuring that private data remains useful to that recipient.

Recently, privacy of this form has been achieved using differential
privacy. Differential privacy 
was originally designed to protect data of individuals in static databases~\cite{Dwork2006_1,Dwork2006_3}. Its goal is to allow accurate statistical analyses of a population while providing strong, provable privacy guarantees to individuals. Differential privacy is appealing because it is immune to post-processing~\cite{Dwork2014algorithmic}, in that post-hoc computations on private data do not weaken privacy's guarantees. For example, filtering private trajectories can be done without harming privacy~\cite{LeNy2014,leny14b}. Differential privacy is also robust to side information~\cite{Kasiviswanathan2008}, in that its privacy guarantees are not defeated by an adversary 
with access to additional information about data-producing entities.  Differential privacy has been extended to dynamical systems~\cite{LeNy2014,le2020differential,rostampour2020privatized} in which trajectory-valued data is protected, and it is this notion of differential privacy that we use.
 
Linear-quadratic (LQ) control is the underlying framework for many existing multi-agent control applications. One example is smart power systems where power forecast, generation, and distribution require access to time series usage data measured by smart meters. In particular, LQ control for load frequency control has been used in power systems~\cite{sargolzaei2016, jiang2012}, with the objective of restoring balance between power consumption and generation. In addition, the work in~\cite{SOUDBAKHSH2017,soudbakhsh2014,singh2016,katewa2019differential} incorporates an LQ control scheme for stability and performance of wide-area power control systems using 
standard phasor measurement units.
Other applications of LQ control include motion planning \cite{chen2019autonomous} and security of cyber physical systems \cite{zhang2018denial}. 
Existing work investigates convergence and performance under various constraints; however, despite the sensitive nature of the data involved, privacy is generally absent in their treatment.

In this paper, we use differential privacy to develop a private multi-agent LQ control framework. Adding privacy noise makes this problem equivalent to a linear quadratic Gaussian (LQG) problem, and the optimal controller will be linear in the expected value of agents' states.
Computing this expected value is a centralized operation, and we therefore augment the network with a cloud computer~\cite{HaleOptimization2014}. In contrast to some existing approaches, the cloud is \emph{not} a trusted third party and does not receive sensitive information from any agent~\cite{hale17}. The cloud instead gathers private information from the agents, estimates their states, and generates optimal inputs. These inputs are transmitted back to the agents, which apply them in their local state updates, and then this process repeats.

\textit{Contributions:} Although there exists a large body of privacy research, privacy parameter interpretation and selection both largely remain the domain of subject matter experts.  Moreover, since offering
privacy guarantees for a control system generally involves sacrificing some level of performance, it is critical to quantify the effects of privacy to rigorously evaluate tradeoffs. Our contributions
are therefore the following: 
\begin{enumerate}[labelindent=0em,labelsep=0.2cm,leftmargin=*]

    \item Developing an algorithm for multi-agent differentially private LQ control. (\textit{Sec.~\ref{sec:results}})
    
	\item  Quantifying sensitive information revealed by bounding filter accuracy in terms of privacy parameters (\textit{Sec.~\ref{sec:effects_of_privacy}})
	\item  Providing quantitative criteria for privacy calibration to trade off information shared and control cost (\textit{Sec. \ref{sec:guidelines}})
	\item  Quantifying the relationship between agents' cost and their privacy levels (\textit{Sec. \ref{sec:costAnalyses}})
\end{enumerate} 

Preliminary versions of this work appeared in~\cite{hale2018privacy,yazdani2020error}. This paper differs from~\cite{hale2018privacy} because it does not rely on a trusted aggregator. 
Further, we quantify the tradeoff between cost and privacy, which was not explored in~\cite{hale2018privacy,yazdani2020error}. 


\textit{Organization:} 
Section~\ref{sec:privacy_review} reviews privacy background.
Section~\ref{sec:problem} defines the private LQG problem, and
Section~\ref{sec:results} solves it. 
Section~\ref{sec:effects_of_privacy} bounds filter error under privacy, and
in Section~\ref{sec:guidelines} we provide guidelines for calibrating privacy. 
Section~\ref{sec:costAnalyses} quantifies the cost of privacy. Next, we provide 
simulations in Section \ref{sec:case}, and then Section~\ref{sec:conclusions} 
concludes.

\section{Review of Differential Privacy}\label{sec:privacy_review}
Differential privacy is a statistical notion of privacy that masks sensitive data while still enabling accurate analyses of it~\cite{Dwork2014algorithmic}. 
It is appealing because post-processing does not weaken its protections. In particular, filtering private data is permitted. 
Moreover, differential privacy is not weakened even if an adversary knows the privacy mechanism used~\cite{Dwork2006_3,Dwork2014algorithmic}.
We briefly review differential privacy here and refer the reader to~\cite{Dwork2014algorithmic,Dwork2006_3,LeNy2014}
for a thorough introduction. 
 


%

We use the ``input perturbation'' approach to differential privacy, which
means that agents add noise directly to their outputs before sharing them. Thus, agents do not ever share sensitive data. 
Privacy guarantees are likewise provided on an individual basis. Formally, each agent's state trajectory will be made approximately indistinguishable from other nearby state trajectories which that agent individually could have produced.

We use the notation $\left[\ell\right]=\{1,\dots,\ell\}$
for $\ell\in\mathbb{N}$.  We consider trajectories of the form
${
Z=\left(Z\left(1\right),Z\left(2\right),\dots\right),
}$
 where $Z\left(k\right)\in\mathbb{R}^{d}$ and~$\|Z(k)\|_2 < \infty$ for all $k$. 
Denote the set\footnote{This notation
comes from the fact that all such~$Z$ have finite truncations
with finite~$\ell_2$-norm. See~\cite{LeNy2014} for additional discussion.}  
of all such sequences by~$\tilde{\ell}^d_2$.

We consider~$N$ agents, and we denote agent $i$'s state
trajectory by $x_{i} \in \tilde{\ell}_2^{n_i}$ for some~$n_i \in \N$. The $k^{th}$ element of
$x_{i}$ is~$x_{i}\left(k\right)\in\mathbb{R}^{n_{i}}$. 
We define our adjacency relation over~$\tilde{\ell}_{2}^{n_i}$. 
\begin{dfn} \label{def:adj} (\emph{Adjacency for Trajectories}) 
	Fix an adjacency parameter $b_i > 0$ for agent $i$.
	Two trajectories~$v_i, w_i \in \tilde{\ell}^{n_i}_2$ are adjacent
	if~$\|v_i - w_i\|_{\ell_2} \leq b_i$. 
    We write~$\textnormal{Adj}_{b_i}(v_i, w_i) = 1$ if~$v_i, w_i$
    are adjacent, and~$\textnormal{Adj}_{b_i}(v_i, w_i) = 0$ otherwise. 
\end{dfn}

This adjacency
relation requires that every agent's state trajectory be made
approximately indistinguishable from all other state trajectories
not more than distance $b_i$ away. Next, we define differential privacy for dynamic
systems. 
This definition considers outputs of agent~$i$ of dimension~$q_i$ at each point in time.
Output signals are in the set~$\tilde{\ell}^{q_i}_2$, over which we
use the~$\sigma$-algebra~$\Sigma^{q_i}_2$ (see~\cite{LeNy2014} for a formal construction). 

\begin{dfn}\label{dfn:differential privacy}
	\emph{(Differential Privacy for Trajectories)} 
	Let~$\epsilon_i>0$ and~$\delta_i\in\left(0,\nicefrac{1}{2}\right)$ be given. A mechanism
	$M$ is $\left(\epsilon_i,\delta_i\right)$-differentially private if,
	for all adjacent ${ x_i,x_i'\in\tilde{\ell}_{2}^{n_i} }$,
	we have
	\begin{equation}
	\mathbb{P}\left[M\left(x_i\right)\in S\right]\le e^{\epsilon_i}\mathbb{P}\left[M\left(x_i'\right)\in S\right]+\delta_i\text{ for all }S\in\Sigma_{2}^{q_i}.
	\end{equation}
\end{dfn}

We enforce this definition with the Gaussian mechanism, defined next. 
We use~$s_1(\cdot)$ for the largest singular value of a matrix,
and~$\mathcal{Q}$ to denote the Gaussian tail integral~\cite{shynk2012probability}.

\begin{lemma}[\emph{Gaussian mechanism}; \cite{LeNy2014}] \label{lem:gaussian mechanism}
	Let agent~$i$ use privacy parameters
	$\epsilon_i>0$ and $\delta_i\in\left(0,\nicefrac{1}{2}\right)$ and adjacency parameter~$b_i > 0$. 
	For outputs~$y_i(k) = C_ix_i(k)$,
	the Gaussian mechanism sets
	$\tilde{y}_i(k)=y_i(k)+v_i(k)$, with~$v_i(k)\sim \mathcal{N}\left(0,\sigma_i^{2}I_{q_i}\right),$
	where $I_{q_i}$ is the $q_i\times q_i$ identity matrix, and  
	$\sigma_i\ge\frac{s_1(C_i)b_i}{2\epsilon_i}(K_{\delta_i}+\sqrt{K_{\delta_i}^{2}+2\epsilon_i}),\text{ with }K_{\delta_i}:=\mathcal{Q}^{-1}\left(\delta_i\right)$. 
	This is~$\left(\epsilon_i,\delta_i\right)$-differentially private with respect to~$\textnormal{Adj}_{b_i}$.
\end{lemma}

For convenience, we set~$\kappa(\delta_i, \epsilon_i) = \frac{1}{2\epsilon_i}(K_{\delta_i}+\sqrt{K_{\delta_i}^{2}+2\epsilon_i})$. 
We use the Gaussian mechanism for the rest of the paper.

\section{Problem Formulation}
\label{sec:problem}
We next introduce the private multi-agent LQG problem. 
Below, we write~$\diag(P_1, \ldots, P_n) := \bigoplus_{i=1}^{n} P_i$
for matrices~$P_1$ through~$P_n$. 

\subsection{Multi-Agent LQ Formulation}
Consider~$N$ agents indexed over $i\in\left[N\right]$.
At time $k$, agent~$i$ has state~$x_{i}\left(k\right)\in\mathbb{R}^{n_{i}}$, with dynamics
\begin{equation}
x_i(k+1) = A_ix_i(k) + B_iu_i(k) + w_i(k),
\end{equation}
where $u_{i}\left(k\right)\in\mathbb{R}^{m_{i}}$,
 $w_{i}\left(k\right)\in\mathbb{R}^{n_{i}}$,
$A_{i}\in\mathbb{R}^{n_{i}\times n_{i}}$, 
and $B_{i}\in\mathbb{R}^{n_{i}\times m_{i}}$.
The distribution of process noise is~${w_{i}\left(k\right)\sim\mathcal{N}\left(0,W_{i}\right)}$, 
where $W_{i}\in\mathbb{R}^{n_{i}\times n_{i}}$ is symmetric and positive definite. 
All process noise terms are independent. 

We define the state~${x(k) = (x_1^T(k)\dots x_N^T(k) )^T \in \mathbb{R}^n}$ and 
control~${u(k) = (u_1^T(k)\dots u_N^T(k) )^T \in \mathbb{R}^m}$, where the dimensions ${n = \sum_{i \in [N]} n_i}$ and ${m = \sum_{i \in [N]} m_i}$. 
Along with
$w(k) = (w_1^T(k), \dots , w_N^T(k))^T \in \R^{n}$, and the matrices
${A = \diag(A_1, \ldots, A_N) \in \R^{n \times n}}$ and ${B = \diag(B_1, \ldots, B_N) \in \R^{n \times m}}$, we have the dynamics
\begin{equation} \label{eq:dynamics}
x(k+1) = Ax(k) + Bu(k) + w(k).
\end{equation}
We consider infinite-horizon problems with cost
\vspace{-10pt}
\begin{multline} \label{eq:cost}
J\left(x,u\right) \!=\!\!\! \lim_{K_{f}\rightarrow\infty} \frac{1}{K_{f}}\mathbb{E}\Bigg\{\!\sum_{k=1}^{K_{f}}\!\left(x\left(k\right) \! - \! \bar{x}\left(k\right)\right)^{T}Q \left(x\left(k\right) \!- \! \bar{x}\left(k\right)\right) \\
+ u\left(k\right)^{T}\!Ru\left(k\right)\!\Bigg\},
\end{multline}
 where $Q \in \R^{n \times n}$ and $R \in \R^{m \times m}$. The vector ${\bar{x}_{i}(k)\in\mathbb{R}^{n_i}}$ is agent $i$'s desired state at time $k$, and we define ${\bar{x}\left(k\right)=(\bar{x}_{1}^T\left(k\right),\dots,\bar{x}_{N}^T\left(k\right))^{T}}$. 
 We make the standard assumption that~$\lim_{k \to \infty} x(k) = \bar{x}$ exists~\cite{Bertsekas2005}.
 
  
 \begin{assumption} \label{as:matrices}
	In the cost $J$, $Q=Q^T \succ 0$ and $R=R^T \succ 0$. The pair $(A, B)$ is controllable,
	and there exists $\Omega$ such that~$Q = \Omega^T\Omega$ and such that the pair $(A, \Omega)$ is observable.
\end{assumption}
 
  Assumption~\ref{as:matrices} is standard in LQ control~\cite{anderson1990,Bertsekas2005,sontag2013}, and it guarantees the existence of a solution to an algebraic Riccati equation that we will encounter below \cite[Chapter 4]{Bertsekas2005}.



\subsection{Differentially Private Information Sharing}\label{subsec:DPcommunications}
The cost $J$ is generally non-separable, which means that it cannot be minimized
by agents using only knowledge of their own states. 
We therefore introduce a cloud computer to aggregate information
and distribute control inputs to the agents. 
The cloud has been used in cyber-physical systems, e.g., in SCADA-based monitoring and state estimation~\cite{SOUDBAKHSH2017,soudbakhsh2014,singh2016,sargolzaei2017}, 
and is a natural choice here. 

 At  time $k$, the cloud requests
from agent $i$ the output value ${y_{i}\left(k\right)=C_{i}x_{i}\left(k\right)}$,
where $C_{i}\in\mathbb{R}^{q_{i}\times n_{i}}$. To protect its state trajectory, agent
$i$ sends a differentially private form of $y_{i}$
to the cloud. 
The cloud uses these private outputs to compute optimal inputs for the agents.
Agents use these inputs in their local state updates, and 
then this process repeats. 

Agent~$i$
adds noise to each $y_{i}\left(k\right)$ before
sending it to the cloud to enforce differential
privacy for~$x_i$. 
Agent~$i$ selects privacy parameters $\epsilon_{i}>0$
and~$\delta_{i}\in\left(0,\nicefrac{1}{2}\right)$ and
adjacency parameter~$b_i > 0$. Then agent~$i$ sends the cloud
\begin{equation}\label{eq:outputt}
\tilde{y}_{i}\left(k\right):=y_{i}\left(k\right)+v_{i}\left(k\right)=C_{i}x_{i}\left(k\right)+v_{i}\left(k\right),
\end{equation}
where the 
privacy\footnote{In Equation~\eqref{eq:outputt}, measurement noise 
inherent to the system
can be included, and all analyses permit this change. 
The form of various Riccati equations will remain the same, with
instances of~$V$ replaced by the sum of~$V$ and the measurement noise covariance matrix. 
However, we focus on bounding the effects of privacy, 
and thus exclude measurement noise.}
noise ${v_{i}\left(k\right)\sim\mathcal{N}\left(0,\sigma^{2}_{i}I_{q_{i}}\right)}$
and ${\sigma_{i}\ge\kappa\left(\delta_{i},\epsilon_{i}\right)s_{1}\left(C_{i}\right)b_{i}}$
from Lemma~\ref{lem:gaussian mechanism}. 
The full privacy vector is
$v\left(k\right)=(v_{1}^T\left(k\right),\dots,v_{N}^T\left(k\right))^{T}$,
and, for all~$k$, we have ${v\left(k\right)\sim\mathcal{N}\left(0,V\right)}$ 
with ${V=\text{diag}\left(\sigma^{2}_{1}I_{q_{1}},\dots,\sigma^{2}_{N}I_{q_N}\right)}$. 
Below we use~${C = \textnormal{diag}(C_1, \ldots, C_N)}$.

Agents' reference trajectories are a source of side information that can reveal their intentions. 
However, agents do not need to reveal their whole reference trajectories to the cloud.
As written, the cost~$J$ depends on~$\bar{x}(k)$ for all~$k$, but, leveraging
the standard average cost-per-stage formulation, one can replace~$\bar{x}(k)$
with~$\bar{x}$ for all~$k$ with no loss of optimality; see~\cite{Bertsekas2005} for a thorough discussion. 
We emphasize that this change is independent of privacy and is a standard approach in infinite-horizon LQG. 
As a result, 
only the limit of agent $i$'s reference trajectory, denoted $\bar{x}_i$, is needed by the cloud to
compute optimal inputs.
Agent $i$ thus privatizes~$\bar{x}_i$ before sharing it\footnote{We expect privatizing the reference limit to be unproblematic in applications in which it only changes the cost incurred, e.g., in applications where states are non-physical quantities. 
	  However, if the reference also encodes some notion of safety that could be affected by privacy, e.g., collision avoidance, the approach we present can be augmented with a low-level reactive controller for that purpose, such as
	  a control barrier function~\cite{ames19}.}.
Agent $i$ selects privacy parameters $\bar{\epsilon}_i > 0$ and $\bar{\delta}_i \in (0,\nicefrac{1}{2})$ 
and adjacency parameter~$\beta_i$. 
Two reference limits~$\bar{x}_i, \bar{x}_i'$ are adjacent if~$\|\bar{x}_i - \bar{x}_i'\|_2 \leq \beta_i$. 
Then privacy noise is added via
$\tilde{x}_i:=\bar{x}_i+\bar{w}_i$. 	
Using the rules for privatizing static data in~\cite[Lemma 1]{LeNy2014},
agent~$i$ generates noise via~$\bar{w}_{i}\sim\mathcal{N}\left(0,\bar{\sigma}^{2}_{i}I_{n_{i}}\right)$, 
with $\bar{\sigma}_{i}\ge\kappa\left(\bar{\delta}_{i},\bar{\epsilon}_{i}\right)\beta_{i}$. 


\color{black}
\begin{problem} \label{prob:newprob}
Let the initial estimate~$\hat{x}(0) = \mathbb{E}[x(0)]$ and
the matrices~$A$,~$B$,~$C$,~$V$, and~$W$ be public information. 
	Minimize
\begin{align}
\tilde{J}\left(x,u\right)=\lim_{K_{f}\rightarrow\infty}\frac{1}{K_{f}}E\Bigg\{\sum_{k=1}^{K_{f}}\left(x\left(k\right)-\right. &\left.\tilde{x}\right)^{T}Q\left(x\left(k\right)-\tilde{x}\right)\\
&+u\left(k\right)^{T}Ru\left(k\right)\Bigg\} 
\end{align}
over all control signals $u$ with $u\left(k\right)\in\mathbb{R}^{m}$, subject to
\begin{align*}
x\left(k+1\right) & =Ax\left(k\right)+Bu\left(k\right)+w\left(k\right)\\
\tilde{y}\left(k\right) & =Cx\left(k\right)+v\left(k\right),
\end{align*}
where agent $i$ has privacy parameters $\left(\epsilon_{i},\delta_{i}\right)$ and $(\bar{\epsilon}_i,\bar{\delta}_i)$.
\end{problem}

%

\section{Private LQG Tracking Control}\label{sec:results}




Problem~\ref{prob:newprob} is an infinite-horizon LQG problem
whose optimal controller~\cite{Yazdani2018} is
\begin{equation}\label{eq:controller}
u^{*}\left(k\right)=L\hat{x}\left(k\right)+Mg,
\end{equation}
where
$M = -\left(R+B^{T}KB\right)^{-1}B^{T}$ and $L = MKA$.
Here, $K$ is the unique positive semidefinite solution to
the discrete algebraic Riccati equation
\begin{equation}\label{eq:K}
K = A^{T}KA-A^{T}KB\left(R+B^{T}KB\right)^{-1}B^{T}KA+Q
\end{equation}
and~$g$ solves
$g=A^{T}[I-KB(R+B^{T}KB)^{-1}B^{T}]g-Q\tilde{x}$.
Without privacy, $g$ would depend on $\bar{x}$, but the cloud only receives its private
form, $\tilde{x}$, and this is what it must use. 

Computing state estimates
for infinite time horizons can use a time-invariant Kalman
filter \cite[Section 5.2]{Bertsekas2005} whose prediction step is 
$\hat{x}^{-}(k+1) = A\hat{x}\left(k\right)+Bu\left(k\right)$.
The \emph{a posteriori} state estimate $ \hat{x}(k) $ is 
computed with
\begin{multline}\label{eq:kalman state update general form}
\!\!\!\!\hat{x}(k \!+ \!1) \!=\! \hat{x}^{-}(k \!+\! 1)
 \!+\! \barsig C^{T}V^{-1}\!\big(\tilde{y}\left(k \!+\! 1\right)-C\hat{x}^{-}(k \!+\! 1)\big),
\end{multline}
where the \emph{a posteriori} error covariance matrix $\barsig$ is given by
$\barsig 
=(C^{T}V^{-1}C+\Sigma^{-1})^{-1}$,
and the \emph{a priori} error covariance $\Sigma$ is the unique positive semidefinite solution
to the discrete algebraic Riccati equation
${\Sigma 
=A(\Sigma^{-1}+ C^{T}V^{-1}C)^{-1}A^{T}+W}$.
The terms $K,L,M,\barsig,\Sigma,$
and $g$ can be all computed beforehand by the cloud to reduce its computational load at
runtime. 

We solve Problem~\ref{prob:newprob} in Algorithm~\ref{alg:main}: 
for all~$i\in[N]$, Algorithm~\ref{alg:main} provides~$\left(\epsilon_i,\delta_i\right)$-differential privacy for agent $i$'s state trajectory and $(\bar{\epsilon}_i,\bar{\delta}_i)$-differential privacy for $\bar{x}_i$.



\begin{algorithm} \label{alg:main}
	\KwData{Public information: $A_i$, $B_i$, $C_i$, $\epsilon_i$, $\delta_i$, $\bar{\epsilon}_i$, $\bar{\delta}_i$, $\hat{x}_i(0)$, $W_i$, and $V_i$ for all~$i$, and $Q$, $R$}
	For all $i$, agent $i$ chooses $(\epsilon_i, \delta_i)$ and
   $(\bar{\epsilon}_i,\bar{\delta}_i)$. It computes $\tilde{x}_i$ and sends it to the cloud \\
	In the cloud, compute~$K$,~$L$,~$M$,~$\Sigma$,~$\barsig$, and~$g$ \\
	\For{$k=0,1,2,\ldots$}{
		\For{$i=1,\ldots,N$}{
			Agent $i$ sends the cloud the private output~$\tilde{y}_i(k) := C_ix_i(k) + v_i(k)$
		}
		In the cloud, compute $u^*(k)$ via \eqref{eq:controller}, send $u^*_i(k)$ to agent $i$ \\
		\For{$i=1,\ldots,N$}{
			Agent $i$ updates its state via $x_i(k+1) = A_ix_i(k) + B_iu^*_i(k) + w_i(k)$
		}
	}
	\caption{Differentially Private LQG (Solution to Problem~\ref{prob:newprob})}
	\label{alg:main}
\end{algorithm}



The feedback control signals $u^*_i$, $i \in [N]$, are computed using estimates of agents' states, and these state estimates are functions of the private output trajectories $\tilde{y}_i$, $i \in [N]$. 
The signals  $ u_i^* $ are thus post-processing on private data and do not reveal agents' state trajectories. 
In addition, knowledge of how $u_i^*$ depends upon the $x_i$'s is equivalent to knowledge of agents' dynamics, which is often assumed to be public information and
is unproblematic for privacy.
Therefore, this use of a feedback controller does not harm privacy.

\section{Quantifying Error Induced by Privacy}\label{sec:effects_of_privacy}
Algorithm~\ref{alg:main} solves Problem~\ref{prob:newprob}, though adding privacy noise
makes it more difficult for the cloud to compute optimal control values. 
Indeed, the purpose
of differential privacy is to protect an agent's state from the cloud, other agents, and any eavesdroppers. 
Thus, the cloud is forced to estimate agents' states to generate control values for
them. Accordingly, 
in this section
we quantify the ability of the cloud to estimate the agents' states as a measure
of the impact of privacy. 

The cloud runs a Kalman filter and computes the input
$u^{*}\left(k\right)$, though privacy noise only affects
the Kalman filter due to the certainty equivalence
principle~\cite{Bertsekas2005}.
We therefore quantify the impact of privacy upon the Kalman
filter in Algorithm~\ref{alg:main} by investigating the best estimate that can be computed with differentially private outputs. 



We proceed by developing trace bounds for the \emph{a priori} error covariance matrix~$\Sigma$ and the \emph{a posteriori} error covariance matrix~$\overline{\Sigma}$, which are, respectively, equal to the 
steady-state mean-square error
(MSE) of the prediction and estimation steps in the Kalman filter. Because the Kalman filter minimizes both of these quantities, lower bounds on them are lower bounds on (asymptotic) MSE for \emph{any} filtering strategy.

We use $ \lambda_{n}(\Upsilon)\le \dots \le \lambda_{1}(\Upsilon) $ to denote the ordered eigenvalues of the matrix $ \Upsilon $.
For simplicity, consider~$C$ diagonal. 
Noting that
${
C^{T}V^{-1}C=\text{diag}\left(\frac{C_{11}^{2}}{\sigma_{1}^{2}},\dots,\frac{C_{nn}^2}{\sigma_{n}^{2}}\right), 
}$
we define
\begin{equation} \label{eq:lu}
l = \argmin\limits_{1\le i\le n}\frac{C_{ii}^{2}}{\sigma_{i}^{2}}, \qquad 
u = \argmax\limits_{1\le i\le n}\frac{C_{ii}^{2}}{\sigma_{i}^{2}}. 
\end{equation}

\begin{theorem}\label{thm:traceboundSigma}
	Suppose every agent	shares its private output trajectory, and the cloud has all public information. Then the steady-state \emph{a priori} MSE of the Kalman filter is bounded via
	\begin{equation} 
	\textnormal{tr}W+\frac{\sigma_{u}^2\textnormal{tr}(A^{T}A)\lambda_{n}(W)}{\sigma_{u}^2+\lambda_{n}(W)C_u^2} \le \textnormal{tr}\Sigma \le \textnormal{tr}W+  \frac{\sigma_{l}^2\textnormal{tr}(A^{T}A)}{C_l^2}
	\end{equation}
and the steady-state \emph{a posteriori} MSE is bounded via
	\begin{equation}
	\frac{n\sigma_{u}^2}{C_u^2+\sigma_{u}^2\lambda_{n}^{-1}(W)} \le \text{tr}\overline{\Sigma}  \leq n\frac{\sigma_{l}^2}{C_l^2},
	\end{equation}
	where~$\sigma_l = \kappa(\delta_l, \epsilon_l)s_1(C_l)b_l$
	and~$\sigma_u = \kappa(\delta_u, \epsilon_u)s_1(C_u)b_u$ are the
	minimum and maximum privacy noise among agents. 
\end{theorem}
\emph{Proof:} See the appendix. \hfill $\blacksquare$

	

These bounds relate privacy to the accuracy of information shared with the cloud and give insight into
differential privacy's protections in conventional estimation-theoretic terms. 
We next leverage these bounds to guide privacy calibration. 

\section{Guidelines for Selecting Privacy Parameters}\label{sec:guidelines}
Calibrating privacy can be challenging. 
The computer science literature has studied this problem~\cite{hsu2014}, though, to the best of our knowledge, there are not control-theoretic guidelines for calibrating 
privacy. Therefore, in this section, we develop such techniques. The privacy parameter $ \delta $ can be interpreted as the probability that~$\epsilon$-differential privacy fails, and is typically~\cite{LeNy2014} chosen in the range $[10^{-5},10^{-1}]$ on this basis. The parameter $ \epsilon $ can be interpreted as the privacy loss of differential privacy, 
and it is typically the parameter to be tuned. We therefore develop guidelines for calibrating~$\epsilon$.

\begin{theorem}\label{thm:guidelineSigmaBar}
	Suppose the cloud has all public information, and agent~$i$ shares its private output trajectory~$\tilde{y}_i$, where~$\tilde{y}_i(k) = C_ix_i(k) + v_i(k)$. 
	Take $ \delta_i \in [10^{-5}, 10^{-1}] $ and set ${\sigma_i = s_1(C_i)\kappa(\delta_i, \epsilon_i)}b_i$. Suppose we want the MSE in the cloud's state estimates to be
	bounded below by~$B_l > 0$ and above by~$B_u > B_l$. 
	These bounds are attained if
	\begin{equation}
	\frac{1}{8}\left(\frac{1+\sqrt{36\eta_{4}+1}}{\eta_{4}}\right)^{2}\le \epsilon_i \le  \frac{1}{\eta_{3}}
	\end{equation}
	for all $i$, where 
	\begin{equation}\label{eq:eta2}
	\eta_{3}\!:=\!\left(\!\frac{B_{l}C_{u}^{2}}{s_1(C_i)^{2}b_i^2\!\left(n\!-\!B_{l}\lambda_{n}^{-1}(W)\right)}\!\right)^{1/2}\!\!\!\!\!\!, \,
	\eta_{4}\!:=\! \left(\!\frac{B_{u}C_{l}^{2}}{ns_1(C_i)^{2}b_i^2}  \right)^{1/2}\!\!\!\!\!\!.
	\end{equation}
\end{theorem}
\emph{Proof:} See the appendix. \hfill $\blacksquare$


Theorem~\ref{thm:guidelineSigmaBar} provides guidelines for choosing~$\epsilon_i$,
which allows agents to make informed decisions for privacy. With this ability, we next examine privacy's impact upon the cost~$J$.

\section{The Control-Theoretic Cost of Privacy}\label{sec:costAnalyses}
Implementing differential privacy adds noise where it would otherwise be absent, and we expect privacy to increase the cost~$J$ relative to a non-private implementation. 
Without any cost considerations, one could add noise of very large variance to provide arbitrarily strong privacy. 
However, private information is used to compute control inputs, which affect future states. 
Thus, there is a need to balance privacy and performance. 
The existing literature has explored several notions of a ``cost of privacy;" LQG minimizes~$J$,
and we therefore compute the increase in~$J$ due to privacy, which offers a ``cost of privacy'' in standard control-theoretic terms. 

\begin{theorem}[\emph{Cost of Privacy}] \label{thm:Cost Diff btwn stochastic and deterministic ref}
	Let $J_{0}(x,u)$ be the cost of Algorithm~\ref{alg:main} without privacy, i.e., with $v_i(k)=\bar{w}_i=0$ for all~$i$ and~$k$. 
	Let~$\tilde{J}(x, u)$ be the cost of Algorithm~\ref{alg:main} with privacy. 
	Then the cost of privacy in LQG, denoted~$\Delta J$, is
	\begin{align}\label{eq:overHeadCost}
	\Delta J\left(x,u\right) &= \tilde{J}(x,u) - J_{0}(x,u)	\\ 
	&=\text{tr}\left(K\Sigma+\left(Q-K\right)\overline{\Sigma}\right) - \text{tr}\left(KW\right) \\
	&+\text{tr}(Q\overline{W})+\text{tr}(H^{T}RH\overline{W}),
	\end{align}
	where~$H = M\big[I - (A + BL)^T\big]^{-1}$. 
\end{theorem}

\emph{Proof:} See the appendix. \hfill $\blacksquare$


	After selecting a privacy level and computing its cost, agents may wish to change their privacy levels to tune costs. 
	For example, agents may choose to relax privacy for significant reductions in~$\Delta J$. 
	A natural way to analyze these changes is with the derivative of the cost of privacy $ \Delta J $ with respect to $ \epsilon_i $; recalling that $ \delta_i $ is typically fixed \emph{a priori}, $ \epsilon_i $ is the parameter to be tuned. For simplicity, we take~$ \bar{\epsilon}_i = \epsilon_{i} = \epsilon $, 
	$ \bar{\delta}_i= \delta_{i} = \delta $,
	and~$s_1(C_i)b_i = \omega$ for all~$i$. 
	
	\begin{theorem}\label{thm:CostRate}
	Let~$A$ be stable. Then the sensitivity of the cost of privacy to changes in privacy is lower-bounded via
	\begin{multline*}
	\frac{d\Delta J}{d\epsilon}\ge\left(-\frac{\omega}{\epsilon}\kappa(\delta,\epsilon)+\frac{\omega}{2\epsilon}\frac{1}{\sqrt{K_{\delta}^{2}+2\epsilon}}\right)\cdot\\
	\Bigg(\!\lambda_{1}(K)\frac{-2\sigma\text{tr}(F^{T}F)}{\lambda_{n}(U)} + 2\sigma\text{tr}Q+2\sigma\text{tr}\left(H^{T}RH\right) + \\
	\text{\ensuremath{\lambda_{1}}}(Q-K)\Bigg[\!\max\! \left\lbrace \frac{-2\sigma\text{tr}(F^{T}F)}{\lambda_{1}(U)}, 0\right\rbrace \lambda_{n}\left(\bar{U}\right)+\text{tr}\left(2\sigma\bar{F}^{T}\bar{F}\right)\Bigg]\Bigg)
	\end{multline*}
	and upper-bounded via
	\begin{multline*}
	\frac{d\Delta J}{d\epsilon}\le \left(-\frac{\omega}{\epsilon}\kappa(\delta,\epsilon)+\frac{\omega}{2\epsilon}\frac{1}{\sqrt{K_{\delta}^{2}+2\epsilon}}\right)\cdot\\\Bigg(\lambda_{n}(K)\max \left\lbrace\frac{-2\sigma\text{tr}(F^{T}F)}{\lambda_{1}(U)}, 0\right\rbrace +2\sigma\text{tr}Q+2\sigma\text{tr}{\left(H^{T}RH\right)} \\ \lambda_{n}(Q-K) \Bigg[\!\frac{-2\sigma\text{tr}(F^{T}F)}{\lambda_{n}(U)}\lambda_{1}\left(\bar{U}\right)+\text{tr}\left(2\sigma\bar{F}^{T}\bar{F}\right)\!\Bigg]\!\Bigg),
	\end{multline*}
	where we use the matrices ${P=C^{T}\left(C\Sigma C^{T}+V\right)^{-1}C\Sigma A^{T}}$, ${U=(A^{T} \!-\! P)(A \!-\! P^{T})-I}$, ${\bar{P}=C^{T}\left(C\Sigma C^{T}+V\right)^{-1}C\Sigma}$, 
	${F = \left(C\Sigma C^{T}+V\right)^{-1}C\Sigma A^{T}}$, ${\bar{F} = \left(C\Sigma C^{T}+V\right)^{-1}C\Sigma}$, and $ \bar{U}=(I-\bar{P})(I-\bar{P}^{T}) $.
	
	\emph{Proof:} See the appendix. \hfill $\blacksquare$

	
	\end{theorem}

		Theorem~\ref{thm:CostRate} explores the continuum of privacy costs that result from varying $\epsilon$, and for a given problem it provides parameter regimes that
		either make it useful or not to relax privacy. 
		One may also wish to enforce hard constraints on performance. 
	Next, we provide guidelines for choosing the privacy parameters $\{\epsilon_{i}\}_{i\in [N]}$ to enforce a desired cost bound.

\begin{theorem}\label{thm:privacyTuningBasedonCost}
	 Suppose a performance requirement is given as a bound on cost by requiring $ \tilde{J} (x,u) \le \alpha $. Take ${\delta_i\in[10^{-5},10^{-1}]}$ and 
	 set $\sigma_i  = s_1(C_i)\kappa(\delta_i, \epsilon_i)b_i$. Then Algorithm~\ref{alg:main} attains $  \tilde{J} (x,u) \le \alpha $ 
	 if, for all~$i$,
	$\epsilon_{i}\ge\frac{1}{8}\left(\frac{1+\sqrt{36\eta_{5}+1}}{\eta_{5}}\right)^{2}$,
	where 
	\resizebox{.5 \textwidth}{!}
	{$ 
	\eta_{5} \!=\! \left[\frac{B_{u}-\lambda_{1}\left(K\right)\textnormal{tr}W-\bar{x}^{T}Q\bar{x}+g^{T}B\left(R+B^{T}KB\right)^{-1}B^{T}g}{s_1(C_i)^2b_i^2\left(\lambda_{1}\left(K\right)\frac{\textnormal{tr}(A^{T}A)}{C_{l}^{2}}+\text{tr}(H^{T}RH)+\text{tr}(Q)\right)}\right]^{\nicefrac{1}{2}}\!\!\!\!\!.$
	}
	\end{theorem}
	\emph{Proof:} See the appendix. \hfill $\blacksquare$

\section{Case Study}\label{sec:case}

Load Frequency Control (LFC)
regulates power flow to different areas while balancing load and generation.
In our framework, each area is an agent, and we consider a system of ten decoupled areas. 
LFC requires transmitting measurements from remote terminal units
(RTUs) to a control center and control signals from the control center
to the plant side. 
This aggregation and communication have well-established privacy
concerns~\cite{McDaniel2000,doe2010}, and we use Algorithm~\ref{alg:main} for it. 
The  continuous time dynamic model of the multi-area LFC system  is given by 
$	\dot{X}(t)=A_{c}X(t)+B_{c}U(t), $
and the matrices $A_c$ and $B_c$ can be found in~\cite{jiang2012}.
The state vector for agent~$i$ is
\[
x_{i}(t)=[\Delta f^{i}(t), \,\, \Delta P_{g}^{i}(t), \,\, \Delta P_{tu}^{i}(t),  \,\, \Lambda^{i}(t)]^{T}\in\mathbb{R}^{4}, 
\]
where $\Delta f^{i}(t),\Delta P_{g}^{i}(t),$ and $\Delta P_{tu}^{i}(t)$
are the frequency deviation, generator
power deviation, and position value of the turbine, respectively. The control input error on the $i$-th power area
is denoted by $\Lambda^{i}(t)=\int_{0}^{t}\vartheta_{i}\Delta f^{i}(s)dt$,
where $\vartheta_{i}$ is the frequency bias factor. 
We simulate~$5$ agents with dynamics of Area~1 from~\cite{jiang2012} and~$5$ agents with dynamics of their Area~2. 

\begin{figure}
	\begin{center}
	\includegraphics[width=.7\columnwidth]{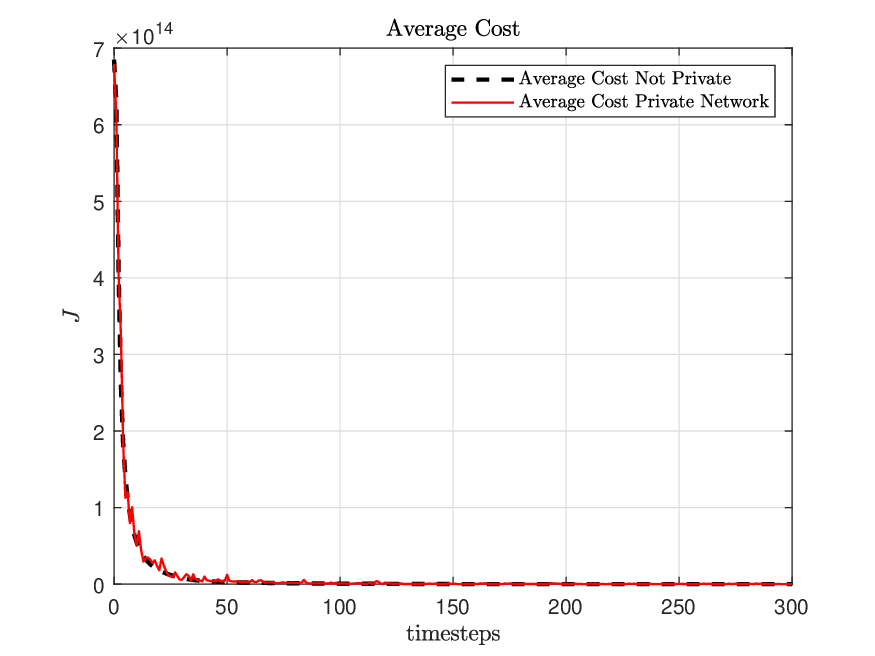}
	\end{center}
	\caption{The time-average cost incurred by ten power system areas with privacy (solid line) and without privacy (dashed line). Privacy increases costs, though these increases become  small relative to the noise-free cost, indicating that privacy's impact is not excessive. This simulation was performed 100 times and 
	the average over these runs is plotted.
	}
	\label{fig:costplot}
\end{figure}

\begin{figure}
\begin{center}		\includegraphics[width=.70\columnwidth]{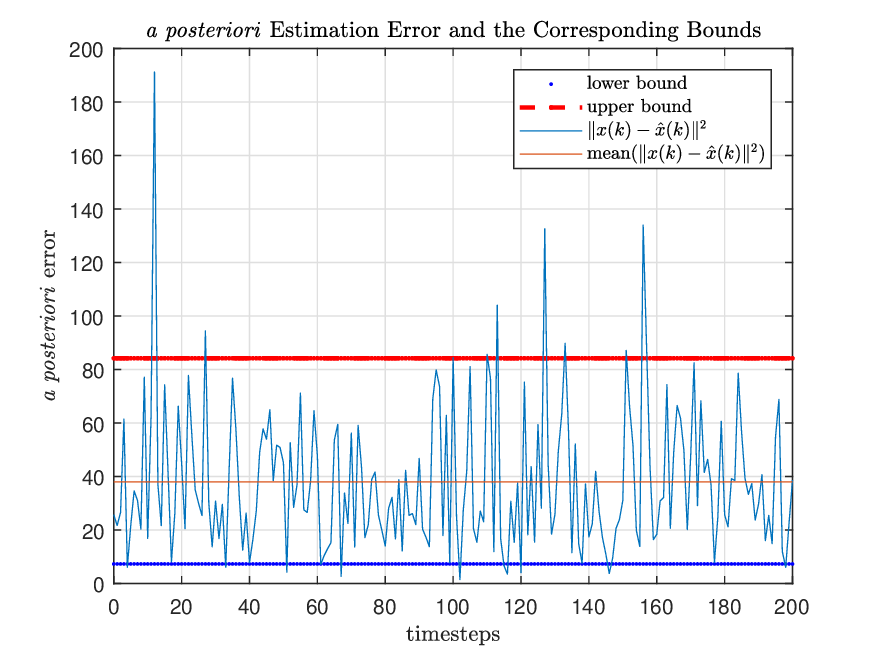}
	\end{center}
	\vspace{0 em}
	\caption{The squared error of the cloud's state estimates (solid line), the lower bound on estimation error in Theorem~\ref{thm:traceboundSigma} (dotted line), 
	and the upper bound on estimation error in Theorem~\ref{thm:traceboundSigma} (dashed line) over $200$ timesteps. We see that instantaneous \emph{a posteriori} error typically obeys our MSE bounds, and on average lies within the bounds.
	}
	\label{fig:post_error}
\end{figure}

We discretize the dynamics of~$X(t)$ with
$A=e^{A_{\mathbf{c}}h}$ and~$B=\int_{0}^{h}e^{A_{\mathbf{c}}\tau}B_{c}d\tau$, 
where $h$ is the sampling period. We have $C=I_{40\times40}$
and $W=I_{40\times40}$. We initialize all states to zero. All areas select identical privacy parameters,
namely, $\left(\epsilon_{i},\delta_{i}\right)=\left(\ln3,0.001\right)$
for all $i$.
In addition, $\Bar{x}_i = 0$ is made private with~$\left(\bar{\epsilon}_{i}, \bar{\delta}_{i}\right)=(\ln 3,0.2)$. 
For the cost, we choose~$Q_{ij} = 100$ for all~$i$ and~$j$, and
we set~$R_{ii} = 100$ and~$R_{ij} = 5$.

\begin{figure}
	\begin{center}
		\includegraphics[width=.7\columnwidth]{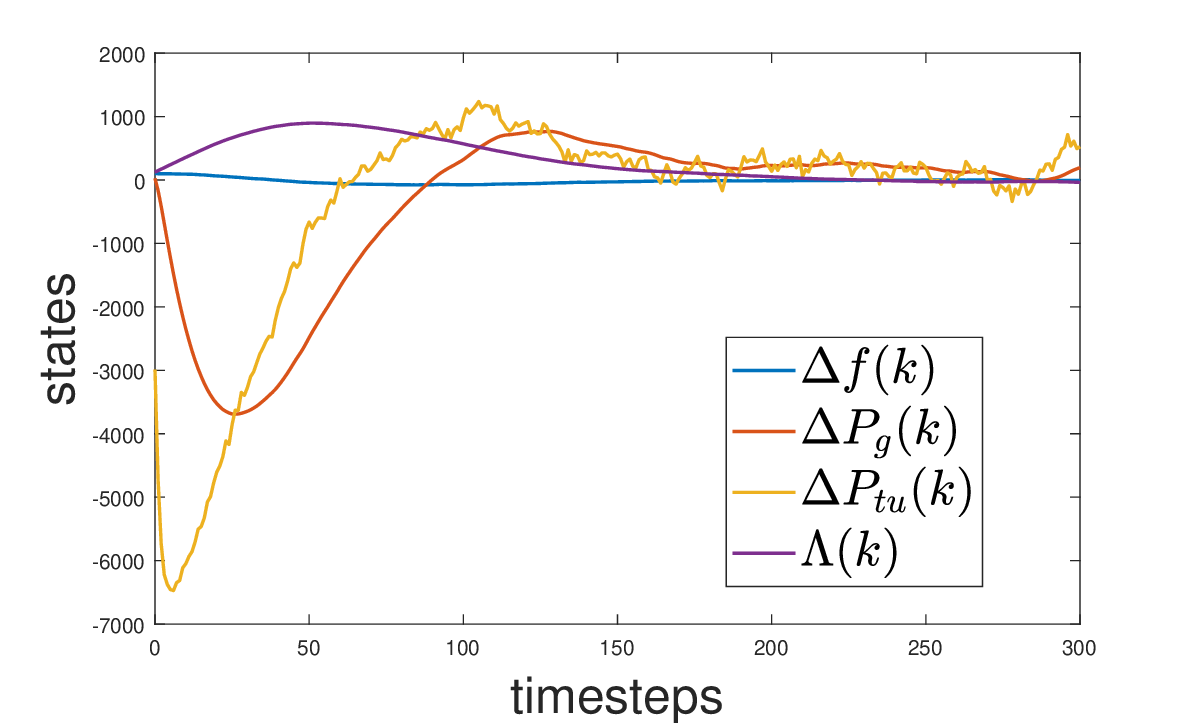}
	\end{center}
	\caption{Even with privacy, control values provided by the cloud are able to regulate an agent's state to remain near its desired trajectory. 
	}
	\label{fig:states}
\end{figure}

The effects of privacy on cost are shown in Figure~\ref{fig:costplot}.
As expected, the cost with privacy is higher than without privacy. However, this increase becomes
relatively modest over time, which indicates that privacy is well-suited to the long-horizon
problems we consider. 
In Figure~\ref{fig:post_error}, we show the instantaneous error of the cloud's state estimates, and we compare that with the bounds in Theorem~\ref{thm:traceboundSigma}; we note that we plot the instantaneous error, but the bounds are for mean-square error. 
As expected, there are ephemeral bound violations by the instantaneous error, and it is shown that  on average, the \emph{a posteriori} error lies within the bounds in Theorem~\ref{thm:traceboundSigma}. This illustrates that privacy
is compatible with the cloud estimating agents' states under privacy. Finally, Figure~\ref{fig:states} illustrates the behaviour of the states of one of the areas.

\section{Conclusions} \label{sec:conclusions}
We have studied distributed linear-quadratic control with differential privacy and bounded the uncertainty and cost induced by privacy. 
Future work will develop differential privacy for other optimal control problems, including in model-free contexts at the intersection
of control and learning~\cite{lewis2012reinforcement}.


\bibliographystyle{IEEEtran}{}
\bibliography{sources}

\section{Appendix}\label{sec:appendix}

The following lemmas will be used in deriving error bounds. 

\begin{lemma}{\cite[Fact 5.12.4]{Bernstein2009}}\label{lem:TraceofProduct}
	Let $\Upsilon$ and $\Theta$ be symmetric $n\times n$ matrices. If $\Upsilon \succ 0$,
then~${\lambda_{n}(\Theta)\textnormal{tr}(\Upsilon)\leq\textnormal{tr}(\Upsilon\Theta)\leq\lambda_{1}(\Theta)\textnormal{tr}(\Upsilon)}$. 
\end{lemma}
\begin{lemma}{\cite[Theorem 8.4.11.]{Bernstein2009}}\label{lem:eigenvalueSplit}
	Let $ \Upsilon $ and $ \Theta $ be $ n\times n $ Hermitian matrices. Then
	\begin{align}\label{key}
	\lambda_1(\Upsilon)+\lambda_n (\Theta)  &\le\lambda_{1}(\Upsilon+\Theta) \le \lambda_1 (\Upsilon) + \lambda_1 (\Theta)\\
	\lambda_n(\Upsilon)+\lambda_n (\Theta)  &\le\lambda_{n}(\Upsilon+\Theta) \le \lambda_ n(\Upsilon) + \lambda_1 (\Theta).
	\end{align}
\end{lemma}

\noindent\textbf{Proof of Theorem~\ref{thm:traceboundSigma}:}
	The steady-state MSE of the Kalman filter's predictions is~$\textnormal{tr}({\Sigma})$.
	Taking the trace of the Riccati equation defining~$\Sigma$, we obtain 
	$\text{tr}\Sigma-\text{tr}W	=\text{tr}\Big[A^{T}A(\Sigma^{-1}+C^{T}V^{-1}C)^{-1}\Big]$,
	where we have used the cyclic permutation property of the trace. Next, we use Lemma~\ref{lem:TraceofProduct} to write
	\begin{align}
	\text{tr}\Sigma &- \text{tr}W \ge\text{tr}(A^{T}A)\lambda_{n}\Big[(\Sigma^{-1}+C^{T}V^{-1}C)^{-1}\Big]\\
	&\ge\frac{\text{tr}(A^{T}A)}{\lambda_{1}(\Sigma^{-1})+\lambda_{1}\big(C^{T}V^{-1}C\big)} =\frac{\text{tr}(A^{T}A)}{\frac{1}{\lambda_{n}(\Sigma)}+\lambda_{1}\big(C^{T}V^{-1}C\big)},
	\end{align}
	where we apply Lemma~\ref{lem:eigenvalueSplit} on the second line to split up the eigenvalues and use the fact that $ \lambda_{1}(\Sigma^{-1})=\nicefrac{1}{\lambda_{n}(\Sigma)} $ in the final step. It is shown in \cite[Theorem 3.1]{Garloff86} that $\Sigma \succeq W$, and therefore $\lambda_n(\Sigma)\ge \lambda_n(W)$. Using this fact and Equation~\eqref{eq:lu} 
	completes the first part of the proof. Similarly, applying Lemmas~\ref{lem:TraceofProduct} and \ref{lem:eigenvalueSplit} to the same Riccati equation, 
	\begin{align}
	\text{tr}\Sigma &- \text{tr}W \le\text{tr}(A^{T}A)\lambda_{1}\Big[(\Sigma^{-1}+C^{T}V^{-1}C)^{-1}\Big]\\
	&\leq \frac{\text{tr}(A^{T}A)}{\lambda_{n}(\Sigma^{-1})+\lambda_{n}(C^{T}V^{-1}C)} \le \frac{\sigma_{l}^2\text{tr}(A^{T}A)}{C_l^2},
	\end{align}
	where the second step uses~$\lambda_{1}(\Upsilon^{-1})=\nicefrac{1}{\lambda_{n}(\Upsilon)}$  and the third step uses Lemma~\ref{lem:eigenvalueSplit} to split the eigenvalues.

	The steady-state MSE of the Kalman filter's state estimates is~$\textnormal{tr}(\overline{\Sigma})$. 
	Using Lemma~\ref{lem:TraceofProduct}, 
	\begin{align*}
	\text{tr}\overline{\Sigma} &\geq \frac{n}{\lambda_{1}(C^{T}V^{-1}C+\Sigma^{-1})} \geq \frac{n}{\lambda_{1}(C^{T}V^{-1}C)+\lambda_{1}(\Sigma^{-1})}\\
	& \geq \frac{n}{\lambda_{1}(C^{T}V^{-1}C)+\lambda_{n}^{-1}(W)} = \frac{n\sigma_{u}^2}{C_u^2+\sigma_{u}^2\lambda_{n}^{-1}(W)},	
	\end{align*}
	where in the second inequality we use Lemma~\ref{lem:eigenvalueSplit} to split the eigenvalues. In the last line,  we use~$\lambda_n(\Sigma)\ge \lambda_n(W)$~\cite[Theorem 3.1]{Garloff86} and Equation~\eqref{eq:lu}. 
	
	Using Lemma~\ref{lem:TraceofProduct}, an upper bound can be derived with
	\begin{align*}
	\text{tr}\overline{\Sigma} &\leq n \lambda_{1}\big((C^{T}V^{-1}C+\Sigma^{-1})^{-1}\big)
	\leq\frac{n}{\lambda_{n}(C^{T}V^{-1}C)},
	\end{align*}
	where we have used Lemma~\ref{lem:eigenvalueSplit} to split the eigenvalues. Using Equation~\eqref{eq:lu} completes the proof. 
\hfill $\blacksquare$

\noindent \textbf{Proof of Theorem~\ref{thm:guidelineSigmaBar}:}
	Choose $ \epsilon_i \ge 	    \frac{1}{8}\left(\frac{1+\sqrt{36\eta_{4}+1}}{\eta_{4}}\right)^{2} $ and solve for $ \eta_{4} $  to get~$\frac{9+\sqrt{2\epsilon_i}}{2\epsilon_i}\le\eta_{4}$.
	Choosing $ \delta_i \in [10^{-5}, 10^{-1}] $ gives ${ K_{\delta_i} \in [1,4.5] }$. Then
	$\frac{2K_{\delta_i}+\sqrt{2\epsilon_i}}{2\epsilon_i}\le \eta_{4}$.
	Because $ \sqrt{\upsilon+\theta} \le \sqrt{\upsilon} + \sqrt{\theta} $, we can lower-bound the left-hand-side to write 
	$\kappa(\delta_i, \epsilon_i) \leq \eta_4$. 
	Squaring, substituting in $ \eta_{4} $, and rearranging we find 
	$s_1(C_i)^2\kappa(\delta_i, \epsilon_i)^2b_i^2 \le\frac{B_{u}C_{l}^{2}}{n}$,
	which implies
	$\sigma_{l}^{2}\le\frac{B_{u}C_{l}^{2}}{n}$. 
	Comparing to Theorem~\ref{thm:traceboundSigma}, we see that~$\textnormal{tr}\overline{\Sigma}\le B_{u}$.

	Next, choose $ \epsilon_i\le  \frac{1}{\eta_{3}} $. Given $ K_{\delta_i}\in [1,4.5] $, we may write
	$
	\eta_{3}\le\frac{K_{\delta_i}}{\epsilon_i}.
	$
	 We substitute for $ \eta_{3} $ and square both sides to write 
	\begin{equation}\label{key}
	\frac{B_{l}C_{u}^{2}}{s_1(C_i)^{2}b_i^2\left(n-B_{l}\lambda_{n}^{-1}(W)\right)}\le \left(\frac{K_{\delta_i}}{\epsilon_i}\right)^{2}. 
	\end{equation}
    Using~$\frac{K_{\delta_i}}{\epsilon_i} \leq \kappa(\delta_i, \epsilon_i)$ and rearranging, we have
	\begin{equation}\label{key}
	\frac{B_{l}C_{u}^{2}}{n-B_{l}\lambda_{n}^{-1}(W)}\le s_1(C_i)^2 \kappa(\delta_i, \epsilon_i)^2 b_i^2. 
	\end{equation}
	This implies
	$\frac{B_{l}C_{u}^{2}}{n-B_{l}\lambda_{n}^{-1}(W)}\le\sigma_{u}^{2}$. 
	Isolating $B_l$ and applying Theorem~\ref{thm:traceboundSigma} implies $\tr \overline{\Sigma} \ge B_l$.
\hfill $\blacksquare$

\noindent \textbf{Proof of Theorem~\ref{thm:Cost Diff btwn stochastic and deterministic ref}:}
	The cost of privatizing~$\bar{x}$ specifically is
	\begin{equation}\label{eq:costofxbar}
	\tr\big([Q + H^TRH]\overline{W}\big),
	\end{equation}
    which is obtained from the authors' technical report in \cite{Yazdani2018} and application of \cite[Equation (318)]{Peterson2012}. Next, using~\cite{Johnson1991}
    and Equation~\eqref{eq:costofxbar},
    the total cost incurred by 
    Algorithm~\ref{alg:main} is
	\begin{align}
	&\tilde{J}(x,u)=J(x,u)+\tr(Q\overline{W}) + \tr (H^T R H\overline{W}) =\\& \lim_{K_{f}\rightarrow\infty}\frac{1}{K_{f}}\sum_{k=1}^{K_{f}}\text{tr}\left(K\Sigma+\left(Q-K\right)\overline{\Sigma}\right)\\&+\lim_{K_{f}\rightarrow\infty}\frac{1}{K_{f}}\sum_{k=0}^{K_{f}-1}\bar{x}^{T}Q\bar{x}-g^{T}B\left(R+B^{T}KB\right)^{-1}B^{T}g \\&+ \tr(Q\overline{W}) + \tr (H^T R H\overline{W}),
	\end{align}
	where the second step follows from \cite[Equation 4.12]{Johnson1991}. Then \eqref{eq:overHeadCost} follows by subtracting the 
	cost without privacy noise.
\hfill $\blacksquare$

	\noindent\textbf{Proof of Theorem~\ref{thm:CostRate}:} 
	Using chain rule we have~$ \frac{d\Delta J}{d\epsilon}=\frac{d\Delta J}{d\sigma}\frac{d\sigma}{d\epsilon} $. For the first term, we have
	\begin{multline}\label{eq:dDeltaJdsigma}
	\frac{d\Delta J(x,u)}{d\sigma}=\text{tr}\left[\frac{d\left(K\Sigma+\left(Q-K\right)\overline{\Sigma}\right)}{d\sigma}\right]\\+\frac{d\left[-\text{tr}\left(KW\right)+\text{tr}(Q\overline{W})+\text{tr}(H^{T}RH\overline{W})\right]}{d\sigma}
	\end{multline}
	where we have used~\cite[Equation 36]{Peterson2012} to move the derivative inside the trace. 
	The matrices~$K$ and~$Q-K$ are symmetric, 
	and~$\text{tr}\frac{d\Sigma}{d\sigma},\text{tr}\frac{d\overline{\Sigma}}{d\sigma}>0$ because filter error monotonically increases with privacy noise.
	Therefore, by Lemma~\ref{lem:TraceofProduct}, the first term in~\eqref{eq:dDeltaJdsigma} can be bounded by
	\begin{multline}\label{eq:dtracedsigma}
	\lambda_{n}(K)\text{tr}\frac{d\Sigma}{d\sigma}+\lambda_{n}(Q-K)\text{tr}\frac{d\overline{\Sigma}}{d\sigma}\le\text{tr}\frac{d\left(K\Sigma+\left(Q-K\right)\overline{\Sigma}\right)}{d\sigma}\\\le\lambda_{1}(K)\text{tr}\frac{d\Sigma}{d\sigma}+\text{\ensuremath{\lambda_{1}}}(Q-K)\text{tr}\frac{d\overline{\Sigma}}{d\sigma}.
	\end{multline}
	By differentiating the discrete algebraic Riccati equation that defines~$\Sigma$, we have
	\begin{multline}\label{key}
	\frac{d\Sigma}{d\sigma}	=A\frac{d\Sigma}{d\sigma}A^{T}-A\frac{d\Sigma}{d\sigma}C^{T}\left(C\Sigma C^{T}+V\right)^{-1}C\Sigma A^{T}\\-A\Sigma C^{T}\left(C\Sigma C^{T}+V\right)^{-1}C\frac{d\Sigma}{d\sigma}A^{T} \\
	+A\Sigma C^{T}\!\!\left(C\Sigma C^{T} \!\!\!+\! V\right)\!^{-1}\!\left(\!C\frac{d\Sigma}{d\sigma}C^{T} \!\!\!+\! 2\sigma I\!\right)\!\!\left(C\Sigma C^{T} \!\!+\! V\right)\!^{-1}C\Sigma A^{T}\!\!.
	\end{multline}
	Taking the trace of both sides and simplifying, we find
	$\text{tr}\left(U\frac{d\Sigma}{d\sigma}\right)=-2\sigma\text{tr}(F^{T}F)$.
	The matrix $ U $ is symmetric, and, because~$A$ is stable, it is positive definite. 
	Therefore by applying Lemma~\ref{lem:TraceofProduct} and simplifying we find 
	\begin{equation}\label{eq:dSigmadsigma}
	\max \left\lbrace \frac{-2\sigma\text{tr}(F^{T}F)}{\lambda_{1}(U)}, 0\right\rbrace 
	\le\text{tr}\frac{d\Sigma}{d\sigma}\le\frac{-2\sigma\text{tr}(F^{T}F)}{\lambda_{n}(U)}.
	\end{equation}
	
	Next we differentiate the equation defining~$\barsig$ to find
	\begin{multline}\label{key}
	\frac{d\overline{\Sigma}}{d\sigma}=\frac{d\Sigma}{d\sigma}-\frac{d\Sigma}{d\sigma}C^{T}\left(C\Sigma C^{T}+V\right)^{-1}C\Sigma\\-\Sigma C^{T}\left(C\Sigma C^{T}+V\right)^{-1}C\frac{d\Sigma}{d\sigma}\\
	+	\Sigma C^{T}\!(C\Sigma C^{T}\!+V)^{-1}\!\left(\!C\frac{d\Sigma}{d\sigma}C^{T} \!\!+\! 2\sigma I\!\right)\!(C\Sigma C^{T}\!+V)^{-1}C\Sigma.
	\end{multline}
	Taking the trace gives
	$\text{tr}\frac{d\overline{\Sigma}}{d\sigma}=\text{tr}\left[\frac{d\Sigma}{d\sigma}\bar{U}\right]+\text{tr}\left(2\sigma\bar{F}^{T}\bar{F}\right)$.
	By substituting the bounds in Equation~\eqref{eq:dSigmadsigma} we get 
	\begin{multline}\label{eq:dSigmabardsigma}
	\max \left\lbrace \frac{-2\sigma\text{tr}(F^{T}F)}{\lambda_{1}(U)}, 0\right\rbrace \lambda_{n}\left(\bar{U}\right)+\text{tr}\left(2\sigma\bar{F}^{T}\bar{F}\right)\le\text{tr}\frac{d\overline{\Sigma}}{d\sigma}\\\le\frac{-2\sigma\text{tr}(F^{T}F)}{\lambda_{n}(U)}\lambda_{1}\left(\bar{U}\right)+\text{tr}\left(2\sigma\bar{F}^{T}\bar{F}\right).
	\end{multline}
	Substituting the results from  Equations~\eqref{eq:dSigmadsigma} and~\eqref{eq:dSigmabardsigma} into~\eqref{eq:dtracedsigma} and assemble the results back in Equation~\eqref{eq:dDeltaJdsigma}, we get
	\begin{multline}\label{eq:dJdsigmaFinal}
	\lambda_{n}(K)\max \left\lbrace \frac{-2\sigma\text{tr}(F^{T}F)}{\lambda_{1}(U)}, 0\right\rbrace +\lambda_{n}(Q-K)\\\Bigg[\frac{-2\sigma\text{tr}(F^{T}F)}{\lambda_{n}(U)}\lambda_{1}\left(\bar{U}\right)+\text{tr}\left(2\sigma\bar{F}^{T}\bar{F}\right)\Bigg]+2\sigma\text{tr}Q+2\sigma\text{tr}{\left(H^{T}RH\right)} \\	\le \frac{d\Delta J}{d\sigma}\le\lambda_{1}(K)\frac{-2\sigma\text{tr}(F^{T}F)}{\lambda_{n}(U)} \\ +
	\lambda_1(Q-K)\Bigg[\max \left\lbrace \frac{-2\sigma\text{tr}(F^{T}F)}{\lambda_{1}(U)}, 0\right\rbrace \lambda_{n}\left(\bar{U}\right)+\text{tr}\left(2\sigma\bar{F}^{T}\bar{F}\right)\Bigg] \\ +2\sigma\text{tr}Q+2\sigma\text{tr\ensuremath{\left(H^{T}RH\right)}}.
	\end{multline}
	
	Next, we 
	observe that $ \frac{d\sigma}{d\epsilon} < 0 $. We multiply Equation~\eqref{eq:dJdsigmaFinal} by $ \frac{d\sigma}{d\epsilon} < 0 $ and this completes the proof.
	\hfill $\blacksquare$

	\noindent\textbf{Proof of Theorem~\ref{thm:privacyTuningBasedonCost}:}
		 Choosing $\epsilon_{i}\ge\frac{1}{8}\left(\frac{1+\sqrt{36\eta_{5}+1}}{\eta_{5}}\right)^{2}$and solving for $\eta_{5}$, we find $ \frac{9+\sqrt{2\epsilon_{i}}}{2\epsilon_{i}}\le\eta_{5} $. Taking $\delta\in[10^{-5},10^{-1}]$ implies $K_{\delta_{i}}\in[1,4.5]$, and as a result we can write ${ \frac{2K_{\delta_{i}}+\sqrt{2\epsilon_{i}}}{2\epsilon_{i}}\le\eta_{5} }$. Using~$\sqrt{\upsilon+\theta}\le\sqrt{\upsilon}+\sqrt{\theta}$ we take
		 $K_{\delta_{i}}$ inside the square root, leading to
		$
		\kappa(\delta_i, \epsilon_i) \le \eta_{5}.
		$
		Expanding, this is equivalent to 
		\begin{multline*}
	    \sigma_i^2 \!\le\! \frac{B_{u} \!-\! \lambda_{1}(K)\textnormal{tr}W \!-\! \bar{x}^{T}Q\bar{x} \!+\! g^{T}B\left(R \!+\! B^{T}KB\right)^{-1}\!\!B^{T}\!g}{\lambda_{1}\left(K\right)\frac{\textnormal{tr}(A^{T}A)}{C_{l}^{2}}+\text{tr}(H^{T}RH)+\text{tr}(Q)}. 
		\end{multline*}
		By rearranging terms and using~$\bar{\sigma}_i = \sigma_i=\sigma_l$, we find 
		\begin{multline*}
		\lambda_{1}\left(K\right)\left[\textnormal{tr}W+\frac{\sigma_{l}^{2}\textnormal{tr}(A^{T}A)}{C_{l}^{2}}\right]+\bar{x}^{T}Q\bar{x}\\-g^{T}B\left(R+B^{T}KB\right)^{-1}B^{T}g
		+\bar{\sigma}_{i}^{2}\left[\text{tr}(H^{T}RH)+\text{tr}(Q)\right]\le B_{u}.
		\end{multline*}
		 Using $\bar{\sigma}_{i}^{2} \leq \lambda_{1}(\overline{W})$ 
		 and Lemma~\ref{lem:TraceofProduct}, we have $\text{tr}(Q\overline{W})+\text{tr}(H^{T}RH\overline{W})\le\text{tr}(Q)\lambda_{1}\left(\overline{W}\right)+\text{tr}(H^{T}RH)\lambda_{1}\left(\overline{W}\right)$
		and we can write 		
		\begin{align}\label{eq:someinequalityFortuning}
		&\lambda_{1}\left(K\right)\left[\textnormal{tr}W+\frac{\sigma_{l}^{2}\textnormal{tr}(A^{T}A)}{C_{l}^{2}}\right]+\bar{x}^{T}Q\bar{x}\\& \!\!\!\!-g^{T}B\left(R+B^{T}KB\right)^{\!-1}\!\!B^{T}g 
		+\text{tr}(Q\overline{W})+\text{tr}(H^{T}RH\overline{W})\le B_{u}.
		\end{align}
		From \cite[Theorem 3.1]{Garloff86} we know that $ Q - K \preceq 0$ and thus~$\lambda_1(Q-K) \leq 0$. 
		Using this in Equation~\eqref{eq:someinequalityFortuning}, we find
		\begin{multline*}
		\lambda_{1}\left(K\right)\left[\textnormal{tr}W+\frac{\sigma_{l}^{2}\textnormal{tr}(A^{T}A)}{C_{l}^{2}}\right]\\+\lambda_{1}\left(Q-K\right)\left[\frac{n\sigma_{l}^{2}}{C_{l}^{2}+\sigma_{l}^{2}\lambda_{n}^{-1}(W)}\right]
		+\bar{x}^{T}Q\bar{x}\\-g^{T}B\left(R+B^{T}KB\right)^{-1}B^{T}g+\text{tr}(Q\overline{W})+\text{tr}(H^{T}RH\overline{W})\le B_{u}.
		\end{multline*}
		Using Theorem~\ref{thm:traceboundSigma}, we can write 
		\begin{multline*}
		\lambda_{1}\left(K\right)\text{tr}\left(\Sigma\right)+\lambda_{1}\left(Q-K\right)\text{tr}\overline{\Sigma}+\bar{x}^{T}Q\bar{x}\\-g^{T}B\left(R+B^{T}KB\right)^{-1}B^{T}g
		+\text{tr}(Q\overline{W})+\text{tr}(H^{T}RH\overline{W})\le B_{u},
		\end{multline*}
		and therefore 
		\begin{multline*}
		\text{tr}\left(K\Sigma\right)+\text{tr}\left[\left(Q-K\right)\overline{\Sigma}\right]+\bar{x}^{T}Q\bar{x}\\-g^{T}B\left(R+B^{T}KB\right)^{-1}B^{T}g
		+\text{tr}(Q\overline{W})+\text{tr}(H^{T}RH\overline{W})\le B_{u},
		\end{multline*}
		and we find 
		$
		\tilde{J}\left(x,u\right)\le\alpha,
		$
		which completes the proof.
	 \hfill $\blacksquare$

\end{document}